\documentclass[number]{elsarticle}
\usepackage{latexsym, amsthm, amssymb, graphicx, anysize, amsmath, mathrsfs}
\usepackage[english]{babel}
\usepackage{hyperref}

\newcommand{\K}{\mathcal{K}}
\newcommand{\Ko}{\mathscr{K}}
\newcommand{\F}{\mathcal{F}}
\newcommand{\Fl}{\mathscr{F}}

\newcommand{\C}{\mathcal{C}}

\newcommand\PG{\mathrm{PG}}

\newcommand\GF{\mathrm{GF}}

\newcommand\Tr{\mathrm{Tr}}

\def\bew{{\bf Proof. }\ }
\def\evb{\hfill \ensuremath{\Box}\\ }

\newtheorem{lemma}{Lemma}
\newtheorem{theorem}{Theorem}

\title{Partial flocks of the quadratic cone yielding Mathon maximal arcs}
\date{}

\author{F. De Clerck\corref{cor1}}
\ead{fdc@cage.ugent.be} 
\address{ Department of  Mathematics, 
Ghent University, Krijgslaan 281 - S22, B-9000 Gent, BELGIUM}

\author{S. De Winter} 
\ead{sgdwinte@me.com} 
\address{ Department of Mathematics, Michigan Technological University, Fisher Hall, 1400 Towsend Drive, Houghton MI 49931, USA}

\author{T. Maes} 
\ead{tmmaes@cage.ugent.be} 
\address{ Department of Mathematics, 
Ghent University, Krijgslaan 281 - S22, B-9000 Gent, BELGIUM}

\cortext[cor1]{Corresponding author}

\begin{document}
\begin{frontmatter}

\begin{abstract}
In \cite{MR2242860} N.~Hamilton and J.~A.~Thas describe a link between maximal arcs of Mathon type and partial flocks of the quadratic cone.  This link is of a rather algebraic nature.  In this paper we  establish 
a geometric connection between these two structures.  We also define a composition on the flock planes and use this to work out an analogue of the synthetic version of Mathon's Theorem (see \cite{DeClerck-DeWinter-Maes:09}).  Finally, we show how it is possible to construct a maximal arc of Mathon type of degree $2d$, containing a Denniston arc of degree $d$ provided that there is a solution to a certain given system of trace conditions.
\end{abstract}
\end{frontmatter}
\section{Introduction}\label{intro}

A $\{k;d\}$-arc $\K$ in a finite projective plane of order $q$ is a non-empty proper subset of $k$ points such that some line of the plane meets $\K$ in $d$ points, but no line meets $\K$ in more than $d$ points.  For given $q$ and $d$, $k$ can never exceed $q(d-1)+d$.  If equality holds $\K$ is called a \textit{maximal arc} of degree $d$, a degree $d$ maximal arc or simply, a maximal $d$-arc.  Equivalently, a maximal arc can be defined as a non-empty proper subset of points such that every line meets the set in $0$ or $d$ points, for some $d$.  The set of points of an affine subplane of order $d$ of a projective plane of order $d$ is a trivial example of a $\{d^2;d\}$-arc. A single point is a $\{1;1\}$-arc of the  projective plane.  We will neglect  these two trivial examples for the rest of this paper. \par
If $\K$ is a maximal $d$-arc 
in a projective plane of order $q$, the set of lines external to $\K$ is a maximal $q/d$-arc 
in the dual plane, called the \textit{dual} of $\K$.  It follows that a necessary condition for the existence of a maximal $d$-arc 
in a projective plane of order $q$ is that $d$ divides $q$.  Denniston showed that this necessary condition is sufficient in the Desarguesian projective plane $\PG(2,q)$ of order $q$ when $q$ is even \cite{MR0239991}.  Ball, Blokhuis and Mazzocca showed that no non-trivial maximal arcs exist in a Desarguesian projective plane of odd order \cite{MR1466573}.  \par
In \cite{MR1883870}, Mathon gave a construction method for maximal arcs in Desarguesian projective planes that generalized the previously known construction of Denniston \cite{MR0239991}.  We will begin by describing this construction method of Mathon.  From now on let $q=2^h$.  \\

\subsection{Mathon maximal arcs}

Let $\Tr$ denote the usual absolute trace map from the finite field $\GF(q)$ onto $\GF(2)$.  We represent the points of the Desarguesian projective plane $\PG(2,q)$ as triples $(a,b,c)$ over $\GF(q)$, and the lines as triples $[u,v,w]$ over $\GF(q)$. A point $(a,b,c)$ is incident with a line $\left[u,v,w\right]$ if and only if $au+bv+cw=0$.  For $\alpha, \beta \in \GF(q)$ such that $\Tr(\alpha \beta)=1$ and $\lambda \in \GF(q)^*= \GF(q) \setminus \{0\}$ we define $F_{\alpha,\beta,\lambda}$ to be the conic
\[
F_{\alpha,\beta,\lambda}=\{(x,y,z):\alpha x^2+xy+\beta y^2+\lambda z^2=0\}.
\]  
Let $\F$ be the set of all such conics.  Remark that all the conics in $\F$ have the point $n(0,0,1)$ as their nucleus and that, due to the trace condition, the line $z=0$ is external to all conics.  \par
For given $\lambda \neq \lambda'$, define a composition
\[
F_{\alpha, \beta, \lambda} \oplus F_{\alpha', \beta', \lambda'}=F_{\alpha \oplus \alpha', \beta \oplus \beta', \lambda \oplus \lambda'}
\]
where the operator $\oplus$ is defined  as follows:
\[
\alpha \oplus \alpha'= \frac{\alpha \lambda + \alpha'\lambda'}{\lambda + \lambda'}, \>\> \beta \oplus \beta'= \frac{\beta\lambda + \beta'\lambda'}{\lambda + \lambda'}, \>\> \lambda \oplus \lambda'= \lambda + \lambda'.
\]
The following lemma was proved by Mathon in \cite{MR1883870}.

\begin{lemma}  \label{discon}
Two non-degenerate conics $F_{\alpha,\beta,\lambda}$, $F_{\alpha',\beta',\lambda'}$, $\lambda \neq \lambda'$ and their composition $F_{\alpha,\beta,\lambda} \oplus F_{\alpha',\beta',\lambda'}$ are mutually disjoint if $\Tr((\alpha \oplus \alpha')(\beta \oplus \beta'))=1$.
\end{lemma}

 Given some subset $\C$ of $\F$, we say $\C$ is \textit{closed} if for every $F_{\alpha,\beta,\lambda} \neq F_{\alpha',\beta',\lambda'} \in \C$, $F_{\alpha \oplus \alpha',\beta \oplus \beta',\lambda \oplus\lambda'} \in \C$.  We can now state Mathon's theorem.

\begin{theorem} [\cite{MR1883870}]  \label{mathon}
Let $\C\subset\F$ be a closed set of conics  in $\PG(2,q)$, $q$ even.  Then the union of the points on the conics of $\C$ together with their common nucleus is a degree $\vert \C \vert +1$ maximal arc in $\PG(2,q)$.
\end{theorem}

We will sometimes call $n$ the nucleus of the maximal arc; by \cite{MR1997407} this is well defined. 
Note that a maximal arc of degree $d$ of Mathon type contains Mathon sub-arcs of degree $d'$ for all $d'$ dividing $d$ (see \cite{MR1883870}). Every maximal arc isomorphic to one as constructed above will be called a {\it maximal arc of Mathon type}. As we mentioned above, Mathon's construction is a generalization of a previously known construction of Denniston.  This can be seen as follows.  Choose $\alpha \in \GF(q)$ such that $\Tr(\alpha)=1$.  Let $A$ be a subset of $\GF(q)^*$ such that $A \cup \{0\}$ is closed under addition.  Then the point set of the conics 
\[
\K_A=\{F_{\alpha,1,\lambda}:\lambda \in A\}
\] 
together with the nucleus $n(0,0,1)$ is the set of points of maximal arc of a degree $\vert A \vert +1$  in $\PG(2,q)$. Every maximal arc isomorphic to such an arc will be called a {\it maximal arc of Denniston type}.  The conics in $\K_A$ are a subset of the \textit{standard pencil of conics} given by
\[
\{F_{\alpha,1,\lambda}:\lambda \in \GF(q)^*\}.
\]  
This pencil partitions the points of the plane, not on the line $z=0$, and distinct from $n(0,0,1)$, 
 into $q-1$ disjoint  conics on the common nucleus $n$. 
   The line $z=0$ is often called the \textit{line at infinity} of the pencil.  It has been proved in 
   \cite{MR1883870} that all degree-4 maximal Mathon arcs are necessarily of Denniston type.  \par

The following lemma was proved by Aguglia, Giuzzi and Korchmaros.  

\begin{lemma}[\cite{MR2443285}]  \label{aggiko}
Given any two disjoint conics $C_1$ and $C_2$ on a common nucleus.  Then there is a unique degree-$4$ maximal arc of Denniston type containing $C_1 \cup C_2$.
\end{lemma}

In \cite{DeClerck-DeWinter-Maes:09} this was generalized to a synthetic version of Mathon's construction. 

\begin{theorem}[Synthetic version of Mathon's theorem]  \label{synthmathon}
Given a degree-$d$ maximal arc $M$ of Mathon type, consisting of $d-1$ conics on a common nucleus $n$, and a conic $C_d$ disjoint from $M$ with the same nucleus $n$, then there is a unique degree-$2d$ maximal arc of Mathon type $\left<M,C_d\right>$ containing $M \cup C_d$.
\end{theorem}

\subsection{Partial flocks}

This subsection will serve as an introduction to partial flocks as well as a brief description of the algebraic link between partial flocks and maximal arcs of Mathon type as it was proved in \cite{MR2242860}.  Suppose that $\Ko$ is a quadratic cone in $\PG(3,q)$ with vertex $x$.  A \textit{partial flock}\index{partial flock} $\Fl$ of $\Ko$ is a set of disjoint (non-singular) conics on the cone $\Ko$.  A partial flock is called \textit{complete}\index{partial flock!complete} if it is not contained in a larger partial flock.  A \textit{flock}\index{flock} $\Fl$ of $\Ko$ is a partial flock of size $q$. The planes containing the conics of the partial 
flock are called the \textit{flock planes}\index{flock planes}.  If all the flock planes  of a partial flock have a line in common, then this partial flock is called \textit{linear}\index{partial flock!linear}.  Flocks are related to some elation generalized quadrangles of order $(q^2,q)$, line spreads of $\PG(3,q)$ and, when $q$ is even, families of ovals in $\PG(2,q)$, called herds (\cite{MR1344552}).  \\

Suppose that the cone $\Ko$ has equation $X_1X_3=X_2^2$.  The vertex is the point $x(1,0,0,0)$ and does not belong to any plane of a (partial) flock $\Fl$.  The conics of $\Fl$ are defined by $k$ planes $V_i$, $i \in \{1, \ldots, k\}$, of which the equations can be written in the form
\begin{align}  \label{GeneralPlanes}
X_0+f(t)X_1+tX_2+g(t)X_3=0,
\end{align}
with $t \in B$, where $B$ is some subset of $\GF(q)$, and $f$ and $g$ are functions from $B$ to $\GF(q)$. 
The property that every two conics of $\Fl$ are disjoint  is equivalent to 
\begin{align}  \label{FlockTrace}
\Tr\Big[ \frac{(f(s)+f(t))(g(s)+g(t))}{(s+t)^2}\Big]=1, \>\>  \forall s,t \in B, s \neq t.
\end{align}

It is well-known (\cite{MR1883870}) that a closed set of conics $\C$, which can be used to construct maximal arcs of Mathon type, may be written in the form
\[
\C=\{(x,y,z): p(\lambda)x^2+xy+r(\lambda)y^2+\lambda z^2=0, \lambda \in A\},
\]
where $A$ is a subset of $\GF(q) \setminus \{0\}$ such that $A \cup \{0\}$ is closed under addition and $p$ and $r$ are functions from $A$ to $\GF(q)$.  \\

N.~Hamilton and J.~A.~Thas proved in \cite{MR2242860} that the functions $p$ and $r$ associated to $\C$ give rise to a partial flock in the following way.  Set $B=A \cup \{0\}$ and define the functions $f$ and $g$ on $B$ by $f(0)=g(0)=0$ and $f(t)=tp(t), g(t)=tr(t)$ for $t \in A$.  Since $A, p$ and $r$ define a closed set of conics we know that 
\begin{align}  \label{ClosedProperty}
\frac{sp(s)+tp(t)}{s+t}=p(s+t) \>\>\>  {\rm and} \>\>\>  \frac{sr(s)+tr(t)}{s+t}=r(s+t),
\end{align} 
for $s,t \in A$, with $s \neq t$.  As $s+t \in A$ the trace condition for the closed set of conics gives us
\begin{eqnarray*}
1=\Tr[p(s+t)r(s+t)]&=&\Tr\Big[ (\frac{sp(s)+tp(t)}{s+t})(\frac{sr(s)+tr(t)}{s+t})\Big]  \\
&=&\Tr\Big[ \frac{(f(s)+f(t))(g(s)+g(t))}{(s+t)^2} \Big] .
\end{eqnarray*}
This implies that $f,g$ and $B$ define a partial flock.  

{}From (\ref{ClosedProperty}) we know that $sp(s)+tp(t)=(s+t)p(s+t)$ and $sr(s)+tr(t)=(s+t)r(s+t)$, or equivalently that $f(s)+f(t)=f(s+t)$ and $g(s)+g(t)=g(s+t)$.  In other words, the functions $f$ and $g$ arising from a closed set of conics are additive on $B$ and also $B$ is closed under addition.  A partial flock with these properties is called an \textit{additive}\index{partial flock!additive} partial flock in \cite{MR2242860}.  

Conversely, suppose an additive partial flock is given with functions $f$ and $g$ on an additive subgroup $B$ of $\GF(q)$.  Now define $A=B \setminus \{0\}$ and functions $p(t)=f(t)/t$ and $r(t)=g(t)/t$, $t \in A$.  It can be checked that these functions have the required trace and closure conditions on $A$ to give a closed set of conics, and hence a maximal arc of Mathon type in $\PG(2,q)$.  Knowing all the above, the following theorem holds.

\begin{theorem} [\cite{MR2242860}]  \label{MathonFlockEquiv}
A degree-$d$ maximal arc of Mathon type gives rise to an additive partial flock of size $d$ of the quadratic cone in $\PG(3,q)$, and conversely.
\end{theorem}

It was also mentioned in \cite{MR2242860} that a partial flock, corresponding to a maximal arc $M$ of degree $d$ of Mathon type is linear if and only if $M$ is of Denniston type.  We will frequently rely on this property throughout this paper.

\section{Projection}  \label{projection}

As became clear in the previous section a maximal arc of degree $d$ of Mathon type gives rise to an additive partial flock of size $d$ of the quadratic cone in $\PG(3,q)$, and conversely.  The link between these two geometric structures is of an algebraic nature and is based on the trace condition of Mathon's construction.  The authors of  \cite{MR2242860} also remark in their paper that a closed set of conics of size $d-1$ on a common nucleus in $\PG(2,q)$, $q$ even, can be projected from a point 
 onto the quadratic cone and in this way induces a partial flock of the quadratic cone.  However, this partial flock does not have as many nice properties as the one arising from the algebraic approach. 

In this section we will establish a more geometric link between the maximal arcs of Mathon type in $\PG(2,q)$ and additive partial flocks in $\PG(3,q)$.  This is done by obtaining a geometric link between the partial flock arising from projection  from a point 
and the additive partial flock. We will see that the relation between the two partial flocks basically is an ``inversion" on the nuclear line of the cone.

Before continuing we first provide a short lemma that guarantees that our projections are well defined.

\begin{lemma}\label{lemma3}
Let $\Ko$ be a quadratic cone with vertex $x$ in $\PG(3,q)$, let $N$ be its nuclear line, and let $\pi$ be any plane not through $x$. Denote $N\cap \pi$ by $n$, and let $p$ be any point on $N$ distinct from $x$ and $n$. Then the projection from $p$ of any conic $C$ in $\pi$ with nucleus $n$ onto the cone $\Ko$ is a conic on $\Ko$. 
\end{lemma}
\bew
First note that every line through $p$ intersects the cone in a unique point. Hence the projection of $C$ results in $q+1$ points on $\Ko$. We need to show they form a conic. 

Consider any plane $\gamma$ in $\PG(3,q)$ not containing $x$ and not containing $p$. Then $\gamma$ clearly intersects $\Ko$ in a conic, and if we project this conic from $p$ onto $\pi$ then we obtain a conic in $\pi$ having $n$ as its nucleus. In this way we obtain $q^2(q-1)$ conics in $\pi$ with nucleus $n$.

On the other hand, in $\PG(2,q)$ every conic with nucleus $(0,0,1)$ is of the form $\alpha X_0^2 +X_0X_1+\beta X_1^2 + \lambda X_2^2=0$ with $\lambda\neq0$ and $\alpha, \beta$ arbitrary elements of $\GF(q)$. Hence there are $q^2(q-1)$ conics having a given point as their nucleus. 

It follows that the conics with nucleus $n$ in $\pi$ are in one-to-one correspondence with the planes not through $x$ or $p$. The lemma follows.
\evb

Now let $M$ be a degree-$d$ maximal arc of Mathon type in the plane $\PG(2,q)$.  Embed $\PG(2,q)$ in $\PG(3,q)$ and assume that $\PG(2,q)$ is the plane with equation $X_0=0$.  To simplify the calculations ahead we will assume that the conics contained in $M$ have equations
\begin{align}  \label{ConicsNew}
\alpha^2X_1^2+X_1X_3+\beta^2X_3^2+\lambda^2X_2^2=0,
\end{align}
with $\alpha, \beta$ and $\lambda$ elements of $\GF(q)$.  Of course the quadratic polynomial $\alpha^2 x^2+x+\beta^2$ has to be irreducible over $\GF(q)$ and this is satisfied if $ \Tr(\alpha^2 \beta^2)=\Tr(\alpha \beta)=1$.  Hence the change of notation does not alter the trace condition.  In the plane $X_0=0$ all conics contained in $M$ have nucleus $n(0,1,0)$.  These conics will sometimes be denoted by $C:(\alpha^2, \beta^2,\lambda^2)$.

Next, let $\Ko$ be a quadratic cone in $\PG(3,q)$.  Suppose the cone $\Ko$ has equation $X_1X_3=X_2^2$.  The vertex is the point $x(1,0,0,0)$ and does not belong to the plane $X_0=0$.  Notice that the conic which is the intersection of $\Ko$ and the plane $X_0=0$ is not contained in $M$ since the elements $\alpha$ and $\beta$ cannot be zero.  It is clear that the nuclear line $N$ is the intersection of the planes $X_1=0$ and $X_3=0$ which is the line with points $(t,0,1,0), t \in \GF(q)$ and the vertex $x$.  Notice that $N$ intersects $X_0=0$ in the point $n(0,0,1,0)$, the common nucleus of all conics in $M$.

Take the point $p(1,0,1,0)$ on the line $N$.  If $\alpha^2X_1^2+X_1X_3+\beta^2X_3^2+\lambda^2X_2^2=0$ is the equation of a conic $C$ in $M$, then, by Lemma \ref{lemma3} it is enough to project three points of the conic  from $p$ on the cone $\Ko$ and to calculate to equation of the plane spanned by the three projected points.  One easily checks that  this  plane has equation 
\begin{align}  \label{Cplanes}
\lambda X_0+\alpha X_1+ (\lambda +1)X_2 + \beta X_3 =0.
\end{align}
These planes will be called \textit{conic planes}. These planes together with the plane $X_0+X_2=0$, called {\em the singular plane},  define a partial flock with $d$ elements  on the cone $\Ko$.  However this partial flock is not additive, opposed to the one defined in   \cite{MR2242860}. Note that in our notation the additive flock is formed by the planes
\begin{align}  \label{HamThasplanes}
X_0+\alpha^2\lambda^2 X_1+ \lambda^2 X_2+ \beta^2\lambda^2 X_3=0
\end{align}
together with the plane $X_0=0$.  

Now, consider the automorphism $\delta \in \mathrm{PGL}(4,q)$ given by
\begin{align*}
X_0 &\rightarrow X_0+X_2,  \\
X_i  &\rightarrow X_i, \>\>  i >0, 
\end{align*}
that fixes the cone $\Ko$.  This automorphism $\delta$ will map the singular plane $X_0 +X_2=0$ on the plane $X_0=0$ while  the conic planes are mapped on 
\begin{align} \label{deltaPlanes}
\lambda X_0+\alpha X_1 +X_2+\beta X_3=0.
\end{align}


The planes found in (\ref{deltaPlanes}) intersect the nuclear line in the points $(1,0,\lambda,0)$.
Next, consider the inversion $\iota$ on the nuclear line  defined by
\begin{align*}
(1,0,y,0) &\mapsto (1,0,1/y,0),  \>\> y\neq 0  \\
(1,0,0,0) &\mapsto (0,0,1,0),  \\
(0,0,1,0) &\mapsto (1,0,0,0). 
\end{align*}
Then $\iota$ induces an involution on the points of the nuclear line, fixing the point $p(1,0,1,0)$.  We can now use $\iota$ to construct a map $\phi$ on each plane $V$ that doesn't intersect any of the points $x(1,0,0,0)$ nor $n(0,0,1,0)$. 
 Each one of these planes intersects the plane $X_0=0$ in a unique line $L$ and the nuclear line in a point $(1,0,y,0)$, $y \neq 0$.  Define the map $\phi$ on the planes $V$ as follows:
\[
\phi(V)=\langle L,\iota(1,0,y,0)\rangle=\langle L,(1,0,1/y,0)\rangle.
\]
Applying $\phi$ on the planes (\ref{deltaPlanes}) results in  the planes 
\begin{align}\label{eq8}
X_0+\alpha\lambda X_1+\lambda X_2+\beta\lambda X_3=0.
\end{align}

Finally, we apply the automorphism $\kappa \in \mathrm{PGL}(4,q)$ given by
\[
(a,b,c,d) \mapsto (a^2,b^2,c^2,d^2),
\]
that fixes the cone $\Ko$.  This yields
\[
X_0+\alpha^2\lambda^2 X_1 +\lambda^2 X_2+\beta^2\lambda^2 X_3=0,
\]
and so, the conic planes found in (\ref{Cplanes}) are mapped on 
 the planes given in (\ref{HamThasplanes}). 

We can summarize as follows.  Let $M$ be  a maximal arc of Mathon type in the plane $X_0=0$ in $\PG(3,q)$ with common nucleus $n(0,1,0)$ and $X_2=0$ as line at infinity, i.e., given the coefficients $(\alpha^2, \beta^2, \lambda^2)$.  Projection from the point $p(1,0,1,0)$ (being a point on the nuclear line) onto the cone $\Ko$ gives rise to a partial flock equivalent to the one with flock planes (\ref{deltaPlanes}) and  $X_0=0$. This partial flock  is not yet additive.  Applying the simple map $\phi$, arising from an inversion on the nuclear line, to these planes, and then the automorphism $\kappa$ gives us the planes (\ref{HamThasplanes}) found in \cite{MR2242860}, i.e., an additive partial flock.  Of course all the above works in both ways. By the way, note that the point $p$ from which we project is completely arbitrarily on the line $N$.  

\section{Plane composition}
It is natural to wonder about the relation between these conic planes and the singular planes and to check whether the equations of these planes can be calculated directly.  We already know from Lemma \ref{aggiko} that, given any two disjoint conics on a common nucleus in a plane, there is a unique third disjoint conic on the same nucleus such that the three conics form a degree-$4$ maximal arc of Denniston type.  This result can be translated to a result concerning conic planes.  \\

We start by introducing a standard equation for  planes not containing the point $p(1,0,1,0)$.  A plane with an equation of the form
\begin{align}  \label{StandEquation}
aX_0+bX_1+(a+1)X_2+cX_3=0,  \>\> a,b,c \in \GF(q)      
\end{align}
is said to have a \textit{standard equation}\index{standard equation}.  This equation is unique when the coefficients of $X_0$ and $X_2$ are distinct, that is when the plane does not contain $p$. 

\begin{lemma}  \label{conicplanescomp} 
Let  $V$ and $W$  be two planes in $\PG(3,q)$, not passing through the vertex $x$ of a cone  $\Ko$, and intersecting the cone $\Ko$  in two disjoint conics. Let $p$ be a point of the nuclear line which is not contained in one of the planes,  then  there is a unique third plane such that the projection of the intersection of these three planes with $\Ko$  from $p$ onto the plane $X_0=0$ induces a degree-$4$ maximal arc of Denniston type.  
\end{lemma}
\bew
Suppose the cone $\Ko$ has equation $X_1X_3=X_2^2$. We may assume that the two conic planes $V$ and $W$ have a standard equation $V: \lambda X_0+\alpha X_1+ (\lambda +1)X_2+ \beta X_3=0$ and $W: \lambda' X_0+\alpha' X_1+ (\lambda' +1)X_2+ \beta' X_3=0$.  These conic planes are associated to the conics $C_1:(\alpha^2, \beta^2,\lambda^2)$ and $C_2:(\alpha'^2, \beta'^2,\lambda'^2)$ in $X_0=0$.  Using Lemma \ref{aggiko} we know that the conic 
\[
C_1 \oplus C_2:\Big(\frac{\alpha^2\lambda^2+\alpha'^2\lambda'^2}{\lambda^2+\lambda'^2}, \frac{\beta^2\lambda^2+\beta'^2\lambda'^2}{\lambda^2+\lambda'^2},\lambda^2+\lambda'^2\Big) 
\]
is the unique conic inducing a degree-$4$ maximal arc of Denniston type containing both $C_1$ and $C_2$.  The unique conic plane corresponding to $C_1 \oplus C_2$ has equation 
\[
V\oplus W:(\lambda +\lambda')X_0+\frac{\alpha\lambda+\alpha'\lambda'}{\lambda+\lambda'} X_1+ (\lambda+\lambda' +1)X_2+ \frac{\beta\lambda+\beta'\lambda'}{\lambda+\lambda'} X_3=0.
\]  
\evb

Notice that the partial flock associated to a maximal arc of Denniston type should be linear.  One easily checks that the three planes in the above lemma indeed have a line in common. Also note that the coefficients  in the standard equation of the plane $V\oplus W$ are obtained using a Mathon composition.  

We know that, if the equation of the conic plane associated to a conic $C:(\alpha^2, \beta^2,\lambda^2)$ is given by $\lambda X_0 + \alpha X_1 +(\lambda+1)X_2+\beta X_3=0$ this equation is standard.  Once the conic plane is set in standard notation we can use the following lemma to determine the singular plane associated to a degree-$4$ maximal arc of Denniston type.

\begin{lemma}  \label{singplanescomp}
Given two conic planes $V$ and $W$ in $\PG(3,q)$, the singular plane inducing the line at infinity of the unique degree-$4$ maximal arc of Denniston type induced by $V$ and $W$ can be found by the sum of the equations of $V$ and $W$.
\end{lemma}
\bew
The conic planes $V: \lambda X_0+\alpha X_1+ (\lambda +1)X_2+ \beta X_3=0$ and $W: \lambda' X_0+\alpha' X_1+ (\lambda' +1)X_2+ \beta' X_3=0$ are associated to the two conics $C_1:(\alpha^2, \beta^2,\lambda^2)$ and $C_2:(\alpha'^2, \beta'^2,\lambda'^2)$ in $X_0=0$.  We are looking for the singular conic in the pencil
\[
\{ \mu C_1+ \nu C_2: \mu, \nu \in \GF(q), \mu, \nu \neq 0 \}.
\]
Since both conics $C_1$ and $C_2$ have $1$ as coefficient of the term $X_1X_3$ the singular conic in the pencil above can be found by simply taking the sum of both conics, i.e., $\mu=\nu=1$.  This gives us $(\alpha^2+\alpha'^2)X_1^2+(\lambda^2+\lambda'^2)X_2^2+(\beta^2+\beta'^2)X_3^2=0$ which is equivalent to 
\begin{align}  \label{Dennistonline}
(\alpha+\alpha')X_1+(\lambda+\lambda')X_2+(\beta+\beta')X_3=0,
\end{align}
yielding the equation of the line at infinity of the unique degree-$4$ maximal arc of Denniston type induced by $V$ and $W$ in the plane $X_0=0$.  Taking the sum of the equations of the two conic planes $V$ and $W$ gives us the plane with equation
\[
(\lambda+ \lambda')X_0+(\alpha+\alpha')X_1+(\lambda+\lambda')X_2+(\beta+\beta')X_3=0.
\]
Intersecting that plane with the plane $X_0=0$ results in the same equation of the line at infinity.  
\evb

Remark that if we take the sum of the conic planes $V$ and $V \oplus W$ in the proof of Lemma \ref{singplanescomp} 
we do not exactly find equation (\ref{Dennistonline}).  However, we do find the same line at infinity.  By multiplying that equation by the right scalar we can always attain equation (\ref{Dennistonline}).
It is clear that we obtain a different singular plane if the equations of the conic planes are not standard.  In that case Lemma \ref{singplanescomp} does not work.  \\

Next we consider the intersections of each of the planes in the partial flock, i.e., the conic planes and the singular plane, with the nuclear line $N$.  We know that $N$ consists of the points $(t,0,1,0), t \in \GF(q)$ and the vertex $x(1,0,0,0)$.  Since the singular plane should always induce a line at infinity on the plane $X_0=0$ in the projection from $p(1,0,1,0)$ we know that this singular plane intersects the nuclear line in the point $p$.  Furthermore, suppose the planes $V, W$ and $V\oplus W$, as seen in the proof of Lemma \ref{conicplanescomp}, are the three conic planes associated to a random degree-$4$ maximal arc of Denniston type.  Their intersections with the nuclear line gives us the points $(\lambda +1,0,\lambda,0)$, $(\lambda' +1,0,\lambda',0)$ and $(\lambda+\lambda' +1,0,\lambda +\lambda',0)$, respectively.  If, to these three points, we add the vertex $x(1,0,0,0)$ we see that, in the $X_2$-component, the elements of the additive group of order $4$ that induce the Denniston $4$-arc are given.   

\section{Analogue of the synthetic theorem}
In Theorem \ref{synthmathon} we described 
a synthetic version of Mathon's theorem.  With the tools given above it is possible to translate this theorem to a theorem concerning partial flocks.  First we will extend the additive linear partial flock of size $4$ corresponding to a degree-$4$ maximal arc of Denniston type.  

\begin{theorem}  \label{SynthpartflockD}
Let $\Fl$ be 
 an additive linear partial flock of size $4$ and let $V'$ be a plane  not containing the point $n(0,0,1,0)$ nor $x(1,0,0,0)$, and such that $V'$ intersects $\Ko$ in a conic disjoint from the elements of $\Fl$.  Then there is a unique additive partial flock of size $8$ containing the conics determined by $V'$ and the four planes defining $\Fl$. 
\end{theorem}
\bew
This follows immediately from the analysis in the previous sections, Theorem \ref{singplanescomp} and  Theorem \ref{synthmathon}.    
\evb

Remark that, if the plane $V'$ in the previous theorem contains the intersection line of the four planes $V_1', \ldots , V_4'$ defining $\Fl$, the partial flock of size $8$ will be linear, and hence will induce a degree-$8$ maximal arc of Denniston type.  

The previous theorem can be generalized to maximal arcs of Mathon type in the following way, the  proof is analogous to the proof of Theorem \ref{SynthpartflockD}..

\begin{theorem}  \label{Synthpartflock}
Let  $\Fl$  be an additive partial flock of size $d$ and let  $V'$ be  a plane not containing the point $n(0,0,1,0)$ such that $V'$ intersects $\Ko$ in a conic disjoint from the elements of $\Fl$.  Then there is a unique additive partial flock of size $2d$ containing the conics determined by $V'$ and the $d$ planes defining $\Fl$. 
\end{theorem}

Using Lemma \ref{singplanescomp} and the equation of the singular planes we can deduce some properties concerning the lines at infinity of a Mathon maximal arc, i.e., the lines at infinity of the Denniston subarcs contained in a Mathon maximal arc.  To simplify the notation we will call these lines \textit{Denniston lines}.

\begin{lemma}  \label{ConcLines1}  
Let $M$ be  a degree-$2d$ maximal  of Mathon type that contains a degree-$d$ maximal arc $D$ of Denniston type.  Then all Denniston lines of $M$ are concurrent. 
\end{lemma}
\bew 
After projection from the point $p(1,0,1,0)$ the maximal arc $D$ gives rise to a linear partial flock on the cone $\Ko$.  In other words, all the planes inducing this partial flock intersect in a common line $L$.  Using Theorem \ref{Synthpartflock} we can choose a suitable plane $V$ to construct the partial flock of size $2d$ that corresponds to the degree-$2d$ maximal arc $M$.  However, this plane $V$ cannot contain the common line $L$ and hence $V$ must intersect $L$ in a point $r$ in $\PG(3,q)$.  Furthermore, using Lemma \ref{singplanescomp}, since all planes defining $M$
 actually are linear combinations of $V$ and the conic planes in the linear partial flock corresponding to $D$, it is clear that $r$ will be contained in all Denniston planes.  Finally, after projection from $p$ on the plane $X_0=0$, we see that all the Denniston lines must be concurrent as they all contain the projection of $r$. 
\evb


Another property regarding the Denniston lines concerns the coefficients $\alpha$ and $\beta$ in the equation of the conics given by (\ref{ConicsNew}). 

\begin{lemma}  \label{AlphaConstant}
The Denniston lines of a maximal arc of Mathon type are concurrent if the coefficient $\alpha$ or $\beta$ is a constant.
\end{lemma}
\bew
Suppose that $\alpha \in \GF(q)$ is a constant in the equation of the conics contained in a maximal arc of Mathon type as given in (\ref{ConicsNew}).  In this case let $V: \lambda X_0+\alpha X_1+ (\lambda +1)X_2+ \beta X_3=0$ and $W: \lambda' X_0+\alpha X_1+ (\lambda' +1)X_2+ \beta' X_3=0$ be two random conic planes.  Using Lemma \ref{singplanescomp} we know that the singular plane induced by $V$ and $W$ has equation
\[
(\lambda +\lambda')X_0+(\lambda +\lambda')X_2+(\beta +\beta')X_3=0.
\] 
It is clear that the point $(0,1,0,0)$ is always contained in this plane.  This implies that all Denniston lines are concurrent.  An analogous argument holds if $\beta$ is a constant.       
\evb

\section{Additive group}  \label{AdditiveGroup}
Consider an additive group $G$ of order $2d$.  In this section we will discuss how, under certain circumstances, it is possible to construct a degree-$2d$ maximal arc $M$ of Mathon type (having $G$ as its related additive group), and containing a degree-$d$ maximal arc of Denniston type.  Let $G:=\{0,1,\lambda_1, \lambda_2, \ldots, \lambda_{2d-2}\}$ and let  $H:=\{0,1,\lambda_1, \ldots, \lambda_{d-2}\}$ be an additive subgroup  of order $d$ of $G$. The elements of $H$ define in the plane $\pi_0$ with equation $X_0=0$ a degree-$d$ maximal arc $D$ of Denniston type  consisting of the conics 
\[
C_{\lambda^2}: X_1^2+X_1X_3+X_3^2+\lambda^2 X_2^2=0,
\]    
with $\lambda^2=1,\lambda_1^2, \ldots, \lambda_{d-2}^2$, together with their common nucleus $n$.  The line at infinity of $D$ is the line $X_0=X_2=0$.  

We choose an element in $G$ that is not contained in $H$, say $\lambda_d$.  It is clear that $H \cup \{\lambda_d\}$ generates $G$.  Because we are trying to construct a degree-$2d$ maximal arc $M$ of Mathon type that contains $D$ we need, using Theorem \ref{synthmathon}, in the plane $\pi_0$ a conic $C: \alpha^2 X_1^2+X_1X_3+\beta^2 X_3^2+\lambda_d^2 X_2^2=0$, with $\alpha, \beta \in \GF(q)$, disjoint from $D$ on the same nucleus $n$.  However, since $M$ contains the degree-$d$ maximal arc $D$ we can assume without loss of generality, using Lemma \ref{ConcLines1} and Lemma \ref{AlphaConstant}, that $\alpha=1$.  It follows that if we can find a suitable element $\beta$, we will be able to construct the entire maximal arc $M$.  

We know that the two conics $C_1$ and $C$ uniquely determine a third conic $C_1 \oplus C$ in order to form a degree-$4$ maximal arc of Denniston type.  This $4$-arc has a unique line at infinity $L$, which is also uniquely determined by $C_1$ and $C$ (see Lemma \ref{singplanescomp}), moreover, $C_1$ and $L$ induce the conic $C$.  This implies that it suffices to determine $L$ in order to find $C$ and thus $M$.  

Since $\alpha=1$ we can assume that $L$ has an equation of the form
\[
\rho X_2+X_3=0,  \quad X_0=0,
\]
$\rho \in \GF(q)$.  The singular plane $S$ associated to $L$ has an equation of the form $AX_0+ \rho X_2 +X_3=0$.  As we know that this plane has to contain the point $p(1,0,1,0)$ 
we find that $A=\rho$.  Hence $S$ has equation
\begin{align}  \label{S1}
\rho X_0 +\rho X_2+X_3=0.
\end{align}
Furthermore, the conic plane that determines $C_1$ has equation $X_0+X_1+X_3=0$ and the conic plane that determines $C$ has equation $\lambda_d X_0+X_1+(\lambda_d +1)X_2 + \beta X_3=0$.  Since these two equations are standard, their sum also provides us with the equation of the associated singular plane $S$.  We find that $S$ must have the equation $(\lambda_d +1)X_0+(\lambda_d +1)X_2+(\beta +1)X_3=0$, or equivalently
\begin{align}  \label{S2}
\frac{\lambda_d +1}{\beta +1}X_0+\frac{\lambda_d +1}{\beta +1}X_2+X_3=0.
\end{align}
{}From (\ref{S1}) and (\ref{S2}) we see that $\rho=(\lambda_d +1)/(\beta +1)$, or equivalently
\begin{align}  \label{BetaFromS}
\beta = \frac{\lambda_d +1}{\rho}+1.
\end{align}

So, the conic $C$ given by the equation 
\[
X_1^2+X_1X_3+\Big(\frac{\lambda_d +1}{\rho}+1\Big)^2 X_3^2+\lambda_d^2 X_2^2=0
\]
 has to be disjoint from the conics $C_{\lambda^2}: X_1^2+X_1X_3+X_3^2+\lambda^2 X_2^2=0$, with $\lambda^2=1,\lambda_1^2, \ldots, \lambda_{d-2}^2$.  Suppose that the point $r(0, x_1,x_2,x_3)$ in $X_0=0$ is a point of $C \cap C_{\lambda^2}$, for some $\lambda \in \{1,\lambda_1^2, \ldots, \lambda_{d-2}^2\}$ .  Then 
 \[
x_2=\frac{(\lambda_d+1)}{\rho(\lambda_d+\lambda)}x_3.
\]
It follows that  the conics $C_{\lambda^2}$ and $C$ are disjoint if and only if
\[
x^2+x+\Big(1+\lambda^2\frac{(\lambda_d+1)^2}{\rho^2(\lambda_d+\lambda)^2}\Big)=0
\] 
has no solutions in $\GF(q)$.  This will be the case if and only if
\begin{align}  \label{TraceRho}
\Tr\Big[ 1+\lambda^2\frac{(\lambda_d+1)^2}{\rho^2(\lambda_d+\lambda)^2} \Big]=1.
\end{align}
Distinguishing the cases $q=2^h$, $h$ odd and $h$ even, we can simplify condition (\ref{TraceRho}) further.  If $h$ is odd we know that $\Tr[1]=1$ and condition (\ref{TraceRho}) is equivalent to
\begin{align}  \label{TraceRhoOdd}  
\Tr\Big[ \frac{\lambda(\lambda_d+1)}{\rho(\lambda_d+\lambda)} \Big]=0. 
\end{align}
On the other hand, if $h$ is even then $\Tr[1]=0$ and we analogously find
\begin{align}  \label{TraceRhoEven}  
\Tr\Big[ \frac{\lambda(\lambda_d+1)}{\rho(\lambda_d+\lambda)} \Big]=1. 
\end{align}
We conclude that all elements $\rho$ that satisfy condition (\ref{TraceRho}) give rise to a suitable element $\beta$ as given in (\ref{BetaFromS}).  Substituting this $\beta$ in the equation $X_1^2+X_1X_3+\beta^2 X_3^2+\lambda_d^2X_2^2=0$, where we assumed $\alpha=1$ as seen above, gives us a conic $C$ which is disjoint from the degree-$d$ maximal arc $D$ and therefore induces a degree-$2d$ maximal arc of Mathon type where the coefficients of the term $X_2^2$ are the squares of the elements in $G \setminus \{0\}$.  \\

Hence, as soon as the above system of trace conditions has a non-trivial solution we can construct a proper maximal degree-$2d$ arc of Mathon type, containing a degree-$d$ maximal arc of Denniston type. In a worst case scenario all the trace conditions could be linearly independent (over $\GF(2)$). In such case, with $q=2^h$ we are guaranteed of the existence of a Mathon maximal arc of degree $2^{\lfloor \log_2(h) \rfloor+1}$ having the prescribed additive group, containing a maximal arc of Denniston type of degree  $2^{\lfloor \log_2(h) \rfloor}$. So in general one should be able to analyze the linear (in)dependence of the trace conditions. Though we do not believe that in general they are all independent, the analysis of dependence seems to be a hard problem, and an interesting topic for future research.

\section*{Acknowledgement} 
Part of this paper has been written while the authors were visiting the Department of Mathematics of Ohio University, Athens, USA.  This research was partly financed by the Research Project  ``Incidence Geometry''  (BOF/GOA/010) at Ghent University. 

\bibliographystyle{plain}

\end{document}